\documentclass[a4paper,oneside,10.5pt]{article}%
\usepackage{amsmath}
\usepackage{amsfonts}
\usepackage{amssymb}
\usepackage{graphicx}%
\providecommand{\U}[1]{\protect \rule{.1in}{.1in}}

\newtheorem{theorem}{Theorem}[section]

\newtheorem{assumption}[theorem]{Assumption}
\newtheorem{example}[theorem]{Example}

\newtheorem{lemma}[theorem]{Lemma}

\newtheorem{remark}[theorem]{Remark}

\numberwithin{equation}{section}

\begin{document}

\title{ Stochastic maximum principle for optimal control problem with
 a stopping time cost functional}
\author{Shuzhen Yang\thanks{ZhongTai Security Institute for Financial Studies, Shandong University,
Jinan, Shandong 250100, PR China. (yangsz@sdu.edu.cn).}
\thanks{This work was supported by the National Natural Science Foundation of China (Grant No.11701330) and Young Scholars Program of Shandong University.} }
\date{}
\maketitle

\textbf{Abstract}: In this study,  we consider an optimal control problem driven by a stochastic differential system with a stopping time terminal cost functional. We establish the stochastic maximum principle for this new kind of an optimal control problem by introducing a discrete terminal system. Finally, we provide an example to describe the main results of this study.

\textbf{Keywords}: stochastic differential equations; stochastic maximum principle; stopping time

{\textbf{MSC2010}: 93E03; 93E20; 60G99

\addcontentsline{toc}{section}{\hspace*{1.8em}Abstract}

\section{Introduction}

In the real world, we usually use the controlled stochastic differential equation (\ref{ine-1}) to describe systems, for example, the optimal investment problem  and   production planning problem in the market.
\begin{equation}
\label{ine-1}
X^u(s)=x+\int_{0}^{s}b(X^u(t),u(t))dt+\int_{0}^{s}\sigma(X^u(t),u(t))dW(t).
\end{equation}
Generally, the cost functional is given as follows:
\begin{equation}
\label{incos-1}
J(u(\cdot))=E\big{[}\displaystyle\int_0^Tf(X^u(t),u(t))dt+\Psi(X^u(T))\big{]}.
\end{equation}

Here, the terminal cost $\Psi(X(\cdot))$ may be dependent on a given stopping time $\tau\leq T$ but not only a constant $T$. In the optimal investment problem, $\tau$ could be used to describe the action of the investor.   Thus, in this study, we will investigate the stochastic maximum principle for the state process (\ref{ine-1})  with the following cost functional:
\begin{equation}
\label{incos-2}
J(u(\cdot))=E\big{[}\displaystyle\int_0^Tf(X^u(t),u(t))dt+\Psi(X^u(\tau))\big{]}.
\end{equation}

Over the past decades, there has been much work concerning the optimality of the state process (\ref{ine-1}) with the cost functional (\ref{incos-1}). Bensoussan \cite{B81} and Bismut \cite{B78} presented a local maximum principle with a convex control set,  while Peng \cite{P90} presented a global  maximum principle with a general control domain that may not be convex. Dynamic programming with related HJB (Hamilton-Jacobi--Bellman) equations and the maximum principle were powerful approaches for solving optimal control problems \cite{YZ99}.
Initially, Zhou \cite{Z95,Z96,Z98}  introduced the concept of near-optimisation. Subsequently, many authors used the near-optimisation method to study the switching LQ-problem, stochastic recursive problem,  linear forward backward stochastic systems, etc. \cite{HLW10a,HHLW11,LYZ05}.

In this study, we find it difficult to investigate the stochastic maximum principle for the state process (\ref{ine-1}) with the cost functional (\ref{incos-2}) by directly applying the spike variation
technique. This is because for a given $\varepsilon>0$ and $E_{\varepsilon}=[v,v+\varepsilon
]\subset [0,T]$,  we cannot verify whether $E_{\varepsilon}\subset [0,\tau]$. Thus, we establish a near-optimal control problem for this study by applying the following multi-time state cost functional:
\begin{equation}
\label{incos-3}
J^n(u(\cdot))=E\big{[}\displaystyle\int_0^Tf(X^u(t),u(t))dt+\sum_{i=1}^{n-1}\Psi(X^u(t_i))1_{\{t_{i-1}\leq \tau < t_{i} \}}+\Psi(X^u(t_n))1_{\{t_{n-1}\leq \tau \leq t_{n} \}}\big{]},
\end{equation}
where $0=t_0< t_1<t_2\cdots<t_n=T$.
Based on the cost functional (\ref{incos-3}), we investigate an optimal control problem for the state process (\ref{ine-1}) under a general control domain. Recently, the optimal control problem under multi-time state cost functional with a convex control domain was studied \cite{Y16c}.

This paper is organised as follows: In Section 2, we present the stochastic optimal
control problem with a stopping time cost functional and investigate a related near-optimal control problem. The stochastic maximum principle for the near-optimal control problem is given in Section 3. In Section 4, we prove the stochastic maximum principle for the optimal control problem under state constraints with stopping time. Subsequently, we provide an example to illustrate the important results of this study. In Section 5, we conclude the optimal control problem and results of this study.

\section{The optimal control problem}

Let $W$ be a $d$-dimensional standard Brownian motion defined on a complete
filtered probability space $(\Omega,\mathcal{F},P,\{ \mathcal{F}(t)\}_{t\geq
0})$, where $\{ \mathcal{F}(t)\}_{t\geq0}$ is the $P$-augmentation of the
natural filtration generated by the Brownian motion $W$. Let ${\displaystyle \tau }$ be the stopping time, which is defined on the filtered probability space $(\Omega,\mathcal{F},P,\{ \mathcal{F}(t)\}_{t\geq 0})$ with values in ${\displaystyle [0,T]}$, where $T>0$ is a given constant, i.e.
$${\displaystyle \{\tau \leq t\}\in {\mathcal {F}}_{t}},\ \ {\displaystyle t\leq T}.$$

Consider the following controlled stochastic differential
equation:
\begin{equation}
d{X}^u(s)=b({X}^u{(s)},u(s))ds+\sigma ({X}^u{(s)},u(s))dW(s), \quad s\in(0,T],\label{ODE_1}%
\end{equation}
with the initial condition $X(0)=x$, where
$u(\cdot)=\{u(s),s\in \lbrack0,T]\}$ is an adaptive control process taking values from set $U$ of $\mathbb{R}^m$. We denote all the adaptive control processes as $ \mathcal{U}[0,T]$, and $b, \sigma$ are the given deterministic functions.

In this study, we consider the following cost functional:
\begin{equation}
J(u(\cdot))=%
E\big{[}{\displaystyle \int \limits_{0}^{T}}
f({X}^u{(t)},u(t))dt+\Psi({X}^u(\tau))\big{]},\label{cost-1}%
\end{equation}
and
\[%
\begin{array}
[c]{l}%
b:\mathbb{R}^m\times U\to \mathbb{R}^m,\\
\sigma:\mathbb{R}^{m}\times U\to \mathbb{R}^{m\times d},\\
f:\mathbb{R}^m\times U\to \mathbb{R},\\
\Psi:\mathbb{R}^{m}\to \mathbb{R}.\\
\end{array}
\]
We set $\sigma=(\sigma^1,\sigma^2,\cdots,\sigma^d)$, and $\sigma^j\in \mathbb{R}^m$ for $j=1,2,\cdots, d$.

Assume $b, \sigma, f$ are uniformly continuous and satisfy the following
 Lipschitz and continuous conditions:

\begin{assumption}
\label{ass-b}There exists a constant $c>0$ such that%
\[%
\begin{array}
[c]{c}%
\left| b(x_{1},u)-b(x_{2},u)\right| +\left| \sigma(x_{1},u)-\sigma(x_{2},u)\right|
 +\left| \Psi(x_{1})-\Psi(x_{2})\right|\leq c\left|x_1-x_2 \right|,\\
\end{array}
\]
$\forall(x_{1},u),(x_{2},u)\in{\mathbb{R}^m}\times U$.
\end{assumption}

\begin{assumption}
\label{assb-b2}
There exists a constant $c>0$ such that
$$
\displaystyle\sup_{u\in\mathcal{U}[0,T]}E\big{[}
\displaystyle\int_0^T\big{[}\left|b(0,u(t))\right|^2+\left|\sigma(0,u(t))\right|^2\big{]}dt  \big{]}\leq c,
$$
where  $\mathcal{U}[0,T]=\{u(\cdot)\in L^2_{\mathcal{F}}(0,T;U)\}.$
\end{assumption}

\begin{assumption}
\label{ass-fai}Let $b, \sigma, f, \Psi$ be twice differentiable at $x$, and their derivatives in $x$ be continuous in $(x,u)$.
\end{assumption}

If Assumptions \ref{ass-b} and \ref{assb-b2} hold, then there exists a unique
solution $X$ for equation (\ref{ODE_1}) \cite{LS78}. By minimising (\ref{cost-1}) over $ \mathcal{U}[0,T]$, any $\bar{u}(\cdot)\in \mathcal{U}[0,T]$
satisfying
\begin{equation}
J(\bar{u}(\cdot))= \underset{u(\cdot)\in\mathcal{U}[0,T]}{\inf}J(u(\cdot)) \label{cost-2}%
\end{equation}
is called an optimal control problem. The corresponding $\bar{u}(\cdot)$ and $\bar{X}(\cdot)$ are called an optimal state trajectory and optimal pair, respectively.

Note that, it is difficult to prove the stochastic maximum principle for the cost functional (\ref{incos-2}) by applying the spike variation technique. For a given $\varepsilon>0$ and $E_{\varepsilon}=[v,v+\varepsilon ]\subset [0,T]$,  we cannot verify whether $E_{\varepsilon}\subset [0,\tau]$. Let $u(\cdot)\in \mathcal{U}[0,T]$ be any given control. Then, we introduce a discrete version of the stopping time $\tau$ as
 $$
 \tau^n=\sum_{i=1}^{n-1}t_i1_{\{t_{i-1}\leq \tau < t_{i} \}}+t_n1_{\{t_{n-1}\leq \tau \leq  t_{n} \}},
 $$
and consider the following multi-time state cost functional:
\begin{equation}
\label{cos-3}
\begin{array}
[c]{rl}
J^n(u(\cdot))=&E\big{[}{\displaystyle \int \limits_{0}^{T}}
f({X}^u{(t)},u(t))dt+\Psi({X}^u(\tau^n))\big{]}\\
=&E\big{[}\displaystyle\int_0^Tf(X^u(t),u(t))dt+\sum_{i=1}^{n-1}\Psi(X^u(t_i))1_{\{t_{i-1}\leq \tau < t_{i} \}}+\Psi(X^u(t_n))1_{\{t_{n-1}\leq \tau \leq t_{n} \}}\big{]},
\end{array}
\end{equation}
where $0=t_0< t_1<t_2\cdots<t_n=T$, and $t_i-t_{i-1}=\frac{T}{n},\ i=1,2,\cdots,n$. When $n$ is large, we have the following lemma:
\begin{lemma}
\label{ner-le1}
Let Assumptions \ref{ass-b} and \ref{assb-b2} hold. Then, for any given $\varepsilon>0$, there exists a large value of $n$ such that
$$
\left|J(\bar{u}(\cdot))-J^n(\bar{u}^n(\cdot))\right|< \varepsilon,
$$
where $(\bar{u}(\cdot),\bar{X}(\cdot))$ is an optimal solution of cost functional $(\ref{cost-1})$, while  $(\bar{u}^n(\cdot),\bar{X}^n(\cdot))$ is an optimal solution of cost functional $(\ref{cos-3})$.
\end{lemma}
\textbf{Proof: }
Based on Assumptions \ref{ass-b} and \ref{assb-b2}, for any given solution $(u(\cdot),X^u(\cdot))$, we have
$$
{E}\left|X^u(t)-X^u(s) \right|\leq C\sqrt{\left| t-s\right|},
$$
where $C$ is a positive constant and $t,s\in [0,T]$. Based on the above inequality, we get
\begin{equation}%
\begin{array}
[c]{rl}
& E\left|X^u(\tau)-\sum_{i=1}^{n-1}X^u(t_i)1_{\{t_{i-1}\leq \tau < t_{i} \}}-X^u(t_n)1_{\{t_{n-1}\leq \tau \leq t_{n} \}}\right|\\
= &E\left|\sum_{i=1}^{n-1}(X^u(\tau)-X^u(t_i))1_{\{t_{i-1}\leq \tau < t_{i} \}}+(X^u(\tau)-X^u(t_n))1_{\{t_{n-1}\leq \tau \leq t_{n} \}}\right|\\
\leq &  C\sqrt{\frac{T}{n}}.
\end{array}
\end{equation}
Notice that, from Assumption \ref{ass-b}, where $\Psi(\cdot)$ is a Lipschitz continuous function, we deduce
 \begin{equation}%
 \label{near-1}
\begin{array}
[c]{rl}
& \left| J(\bar{u}(\cdot))-J^n(\bar{u}(\cdot))  \right|\leq C\sqrt{\frac{T}{n}},\\
&  \left| J(\bar{u}^n(\cdot))-J^n(\bar{u}^n(\cdot))  \right|\leq C\sqrt{\frac{T}{n}}.\\
\end{array}
\end{equation}

Given that,  $(\bar{u}(\cdot),\bar{X}(\cdot))$ is an optimal solution of cost functional $(\ref{cost-1})$ and  $(\bar{u}^n(\cdot),\bar{X}^n(\cdot))$ is an optimal solution of cost functional $(\ref{cos-3})$, we get
 \begin{equation}%
 \label{near-2}
\begin{array}
[c]{rl}
&  J(\bar{u}(\cdot))\leq J(\bar{u}^n(\cdot)) , \\
&  J^n(\bar{u}^n(\cdot))\leq J^n(\bar{u}(\cdot)). \\
\end{array}
\end{equation}
Combining equations (\ref{near-1}) and (\ref{near-2}), we obtain
$$
\left| J(\bar{u}(\cdot))-J^n(\bar{u}^n(\cdot))  \right|\leq C\sqrt{\frac{T}{n}},
$$
where $n= [\frac{C^2T}{\varepsilon^2}]+1$. Hence, we complete the proof. $\ \ \ \ \ \ \ \ \Box$
\begin{remark}
\label{rema-1} For a large value of $n$, Lemma \ref{ner-le1} shows that the minimum value of cost functional (\ref{cost-1}) is close to the minimum value of cost functional (\ref{cos-3}). Thus, we call the  state process (\ref{ine-1})  under cost functional (\ref{cos-3}) a near-optimal control problem of the state process (\ref{ine-1}) under cost functional (\ref{cost-1}). In the following section, we will consider the optimal control problem with cost functional (\ref{cos-3}).

\end{remark}

\section{Stochastic maximum principle}

For notation simplicity, we set
  \begin{equation}%
\begin{array}
[c]{rl}
&  \phi(X^u(t_i))=\Psi(X^u(t_i))1_{\{t_{i-1}\leq \tau < t_{i} \}},\ i=1,2,\cdots,n-1,\\
&  \phi(X^u(t_n))=\Psi(X^u(t_n))1_{\{t_{n-1}\leq \tau \leq t_{n} \}}.
\end{array}
\end{equation}
In this section, we will study the stochastic maximum principle for the following cost functional:
\begin{equation}
\label{costm-1}
J(u(\cdot))=E\big{[}\displaystyle\int_0^Tf(X^u(t),u(t))dt+\displaystyle\sum_{j=1}^n\phi(X^u(t_j))\big{]},
\end{equation}
with state equation (\ref{ODE_1}). Notice that, we consider a general control domain $U$ that need not be convex. The main difficulty is to investigate the variational and adjoint equations. We introduce the first-order and second-order adjoint equations as follows:

The first-order adjoint equations are%
\begin{equation}%
\begin{array}
[c]{rl}%
-d{p}(t)&= \{b_x(\bar{X}{(t)},\bar{u}(t))^{\rm{T}}p(t)+ \displaystyle\sum_{j=1}^d\sigma_x(\bar{X}{(t)},\bar{u}(t))^{\rm{T}}q^j(t) \\
               &-f_x(\bar{X}{(t)},\bar{u}(t))\}dt-q(t)dW(t),\ t\in(t_{i-1},t_{i}),\\
p(t_{i}) &=-\phi_x(\bar{X}(t_i))+p(t_{i}^{+}),\text{ \ }i=1,2,\cdots,n,
\end{array}
\label{prinm-1}%
\end{equation}
where ``$\rm{T}$" means the transform of a vector or matrix,  $t_{i}^{+}$ is the right limit of $t_{i}$, and $p(t_{n}^{+})=0$.

We define the following Hamiltonian function:
\begin{equation*}
H(x,u,p,q)=b(x,u)^{\rm{T}}p+\displaystyle\sum_{j=1}^d\sigma^j (x,u)^{\rm{T}}q^j-f(x,u),\text{ \  \ }%
\end{equation*}
where $(x,u,p,q)\in \mathbb{R}^m\times U\times \mathbb{R}^m\times \mathbb{R}^{m\times d}$.

The second-order adjoint equations are
\begin{equation}%
\begin{array}
[c]{ll}%
-dP(t)= & \bigg{\{}b_x(\bar{X}{(t)},\bar{u}(t))^{\rm{T}}P(t)+P(t)b_x(\bar{X}{(t)},\bar{u}(t))\\
              &+ \displaystyle\sum_{j=1}^d\sigma_x^j(\bar{X}{(t)},\bar{u}(t))^{\rm{T}} P(t)\sigma_x^j(\bar{X}{(t)},\bar{u}(t)) \\
               &+\displaystyle\sum_{j=1}^d\big{[}
               \sigma_x^j(\bar{X}{(t)},\bar{u}(t))^{\rm{T}}Q^j(t)+Q^j(t)\sigma_x^j(\bar{X}{(t)},\bar{u}(t))\big{]}\\
               & +H_{xx}(\bar{X}(t),\bar{u}(t),p(t),q(t))\bigg{\}}dt -Q(t)dW(t),\ t\in(t_{i-1},t_{i}), \\
-P(t_i)= & \phi_{xx}(\bar{X}(t_i))-P(t_i^+),
\end{array}
\label{prinm-2}%
\end{equation}
where $P(t_n^+)=0$.

The main result of this section is the following theorem:
\begin{theorem}
\label{Maximumprinciple} Let Assumptions \ref{ass-b}, \ref{assb-b2} and \ref{ass-fai} hold,
and $(\bar{u}(\cdot),\bar{X}(\cdot))$ be an optimal pair of (\ref{costm-1}).
Then, there exists  $(p(\cdot),q(\cdot))$ and $(P(\cdot),Q(\cdot))$ satisfying the series of first-order adjoint equations (\ref{prinm-1}) and second-order adjoint equations (\ref{prinm-2}) such that
\begin{equation}%
\begin{array}
[c]{ll}%
&H(\bar{X}(t),\bar{u}(t),p(t),q(t))-H(\bar{X}(t),u,p(t),q(t)))\\
\geq &\displaystyle\frac{1}{2} \displaystyle\sum_{j=1}^d\big{[}\sigma^j(\bar{X}{(t)},\bar{u}(t))-\sigma(\bar{X}{(t)},u)\big{]}^{\rm{T}}P(t)
 \big{[}\sigma^j(\bar{X}{(t)},\bar{u}(t))-\sigma(\bar{X}{(t)},u)\big{]},
\end{array}
\label{prinm-3}%
\end{equation}
for any $u\in U$ and $t \in(t_{i-1},t_i)$, $i=1,2,\cdots,n$.
\end{theorem}

Let $(\bar{u}(\cdot),\bar{X}(\cdot))$ be the given
optimal pair of cost functional (\ref{costm-1}). Let $\varepsilon>0,$ and $E_{\varepsilon}=[v,v+\varepsilon
]\subset (t_{i-1},t_i)$, for some $i\in\{1,2,\cdots,n\}$. Let $u(\cdot)\in \mathcal{U}[0,T]$ be any given control. By applying the spike
variation technique, we define the following control:%
\[
u^{\varepsilon}(t)=\left \{
\begin{array}
[c]{cl}%
\bar{u}(t), & \text{if }t\in \lbrack0,T]\backslash E_{\varepsilon},\\
u(t), & \text{if }t\in E_{\varepsilon},
\end{array}
\right.
\]
where $u^{\varepsilon}(\cdot)\in \mathcal{U}[0,T]$. The following Lemma is
useful for proving Theorem \ref{Maximumprinciple}.

\begin{lemma}
\label{lem-2} Let Assumptions \ref{ass-b} and \ref{assb-b2} hold, and
$X^{\varepsilon}(\cdot)$ be the solution of equation (\ref{ODE_1}) under the
control $u^{\varepsilon}(\cdot)$, and $y^{}(\cdot)$, $z^{}(\cdot)$  be the solutions
of the following equations:
\begin{equation}
\label{aprom-1}
\begin{array}
[c]{cl}%
d{y}(t)= & b_x(\bar{X}{(t)},\bar{u}(t))y(t)dt+\displaystyle\sum_{j=1}^d\big{[} \sigma_x^j(\bar{X}{(t)},\bar{u}(t))y(t)\\
&+\sigma^j(\bar{X}{(t)},u^{\varepsilon}(t))
-\sigma^j(\bar{X}{(t)},\bar{u}(t))\big{]}dW^j(t), \\
y(0)= & 0,\quad t\in (0,T],
\end{array}
\end{equation}
and
\begin{equation}%
\label{aprom-2}
\begin{array}
[c]{ll}%
d{z}(t)= & \big{[}b_x(\bar{X}{(t)},\bar{u}(t))z(t)
+\frac{1}{2}b_{xx}(\bar{X}{(t)},\bar{u}(t))(y(t))^2\\
&+b(\bar{X}{(t)},u^{\varepsilon}(t))-b(\bar{X}{(t)},\bar{u}(t))\big{]}dt \\
&+\displaystyle\sum_{j=1}^d\big{[}   \sigma_x^j(\bar{X}{(t)},\bar{u}(t))z(t)+\frac{1}{2}\sigma_{xx}^j(\bar{X}{(t)},\bar{u}(t))(y(t))^2\\
&+(\sigma_x^j(\bar{X}{}(t),u^{\varepsilon}(t))-
\sigma_x^j(\bar{X}{(t)},\bar{u}(t)))y(t)\big{]}dW^j(t),\\
z(0)= & 0,\quad t\in (0,T].
\end{array}
\end{equation}
Then,%
\begin{equation}%
\begin{array}
[l]{l}%
\max_{t\in \lbrack0,T]}E\left \vert y^{}(t)\right \vert
=O(\varepsilon^{\frac{1}{2}}),\\
\max_{t\in \lbrack0,T]}E\left \vert z(t)\right \vert =O(\varepsilon),\\

\max_{t\in \lbrack0,T]}E\left \vert X^{\varepsilon}(t)-\bar{X}(t)-y^{
}(t)\right \vert =O(\varepsilon),  \\

\max_{t\in \lbrack0,T]}E\left \vert X^{\varepsilon}(t)-\bar{X}(t)-y^{
}(t)-z(t)\right \vert =o(\varepsilon),
 \end{array}
 \label{varm-1}%
\end{equation}
and%
\begin{equation}%
\begin{array}
[c]{rl}
& J(u^{\varepsilon}(\cdot))-J(\bar{u}(\cdot))\\
= &
\displaystyle \sum \limits_{i=1}^{n}E\left[\phi_x(\bar{X}(t_i))(y(t_{i})+z(t_{i})) +\phi_{xx}(\bar{X}(t_i))y^{}(t_{i})y^{}(t_{i})\right]\\
& +
E{\displaystyle \int \limits_{0}^{T}}
\big{[}f_x(\bar{X}{(t)},\bar{u}(t))(y(t)+z(t))+\frac{1}{2}f_{xx}(\bar{X}{(t)},\bar{u}(t))(y(t))^2\\
&+f(\bar{X}{(t)},u^{\varepsilon}(t))-f(\bar{X}{(t)},\bar{u}(t))\big{]}dt+o(\varepsilon).
\end{array}
\label{varm-3}%
\end{equation}

\end{lemma}
\textbf{Proof: }
Applying the technique in Lemma 1 of \cite{P90}, we can prove equation (\ref{varm-1}).

Note that%
\begin{equation}%
\begin{array}
[c]{rl}
& J(u^{\varepsilon}(\cdot))-J(\bar{u}(\cdot))\\
= &\displaystyle\sum_{i=1}^nE\big{[} \phi(X^{\varepsilon}(t_i))-\phi(\bar{X}(t_i))+{\displaystyle \int \limits_{0}^{T}}
\big{(}f(X{}^{\varepsilon}(t),u^{\varepsilon}(t))-f(\bar{X}{(t)},\bar{u}(t))\big{)}dt\big{]}.\\
\end{array}
\label{valuem-0}%
\end{equation}
Applying equation (\ref{varm-1}), we get
\begin{equation}%
\begin{array}
[c]{rl}
& J(u^{\varepsilon}(t))-J(\bar{u}(t))\\
= &E\big{[} \Psi_x(\bar{X}(t_n))(y(t_n)+z(t_n))+
\displaystyle\sum_{i=1}^{n}(y(t_i)+z(t_i))+\frac{1}{2}\Psi_{xx}(\bar{X}(t_n))(y(t_n),y(t_n))\\
&+{\displaystyle \int \limits_{0}^{T}}
\big{(}f_x(\bar{X}{(t)},\bar{u}(t))(y(t)+z(t))+\frac{1}{2}f_{xx}(\bar{X}{(t)},\bar{u}(t))(y(t))^2\\
&+f(\bar{X}{(t)},u^{\varepsilon}(t))-f(\bar{X}{(t)},\bar{u}(t))\big{)}dt\big{]}+o(\varepsilon).\\
\end{array}
\label{valuem-1}%
\end{equation}
This completes the proof. $\ \ \ \ \ \ \ \ \Box$
\bigskip

Based on the above Lemma, we now carry out the proof for Theorem
\ref{Maximumprinciple}.

\textbf{Proof of Theorem} \ref{Maximumprinciple}. For$\ t\in (
t_{i-1},t_{i}),$ applying the differential chain rule to $p(t)^{\rm{T}}(y^{
}(t)+z(t))$, and by Assumption \ref{ass-fai}, we have%
\begin{equation}%
\begin{array}
[c]{rl}
&E\big{[} p(t_{i})^{\rm{T}}(y^{}(t_{i})+z(t_i))-p(t_{i-1}^{+})^{\rm{T}}(y^{}(t_{i-1})+z(t_{i-1})\big{]}\\
= & E\big{[}-\phi_x(X(t_i))(y^{}(t_{i})+z(t_{i}))+p(t_{i}^{+})(y^{}(t_{i})+z(t_{i}))
-p(t_{i-1}^{+})(y^{}(t_{i-1}^{})+z(t_{i-1}^{}))\big{]}\\
= &
E{\displaystyle \int \limits_{t_{i-1}}^{t_{i}}}
\bigg{[}f_x(\bar{X}{(t)},\bar{u}(t))(y_{}
^{}(t)+z(t)) \\
&+\displaystyle\frac{1}{2}p(t)^{\rm{T}}b_{xx}(\bar{X}(t),\bar{u}(t))(y(t))^2
+p(t)^{\rm{T}}(b(\bar{X}{(t)},u^{\varepsilon}(t))-b(\bar{X}{(t)},\bar{u}(t)))
\\
&+\displaystyle\sum_{j=1}^d\big{(}\frac{1}{2}q^j(t)^{\rm{T}}\sigma_{xx}^j(\bar{X}(t),\bar{u}(t))(y(t))^2
+q^j(t)^{\rm{T}}(\sigma(\bar{X}{(t)},u^{\varepsilon}(t))-\sigma(\bar{X}{(t)},\bar{u}(t)))\big{)}\bigg{]}dt.
\end{array}
\label{maxm-1}%
\end{equation}
Adding $i$ on both sides of equation (\ref{maxm-1}), we have
\[%
\begin{array}
[c]{rl}
&
E{\displaystyle \sum \limits_{i=1}^{n}}
\big{[}(p(t_{i})(y^{}(t_{i})+z(t_i))-p(t_{i-1}^{+})(y^{}(t_{i-1})+z(t_{i-1}))\big{]} \\
= &
E{\displaystyle \sum \limits_{i=1}^{n}}
\big{[}-\phi_x(\bar{X}(t_i))(y(t_{i})+z(t_{i}))\\
&+p(t_{i}^{+})(y^{}(t_{i})+z(t_{i}))
-p(t_{i+1}^{+})(y^{}(t_{i+1})+z(t_{i+1}))\big{]} \\
= &
E\big{[}-\displaystyle \sum \limits_{i=1}^{n}
\phi_x(\bar{X}(t_i))(y^{}(t_{i})+z(t_i))\big{]}\\
= &
E{\displaystyle \sum \limits_{i=1}^{n}\displaystyle \int \limits_{t_{i-1}}^{t_{i}}}
\bigg{[}f_x(\bar{X}{(t)},\bar{u}(t))(y_{}
^{}(t)+z(t)) \\
&+\frac{1}{2}p(t)^{\rm{T}}b_{xx}(\bar{X}(t),\bar{u}(t))(y(t))^2
+p(t)^{\rm{T}}(b(\bar{X}{(t)},u^{\varepsilon}(t))-b(\bar{X}{(t)},\bar{u}(t)))
\\
&+\displaystyle\sum_{j=1}^d\big{(}\frac{1}{2}q^j(t)^{\rm{T}}\sigma_{xx}^j(\bar{X}(t),\bar{u}(t))(y(t))^2
+q^j(t)^{\rm{T}}(\sigma(\bar{X}{(t)},u^{\varepsilon}(t))-\sigma(\bar{X}{(t)},\bar{u}(t)))\big{)}\bigg{]}dt\\
= &
E{\displaystyle \int \limits_{0}^{T}}
\bigg{[}f_x(\bar{X}{(t)},\bar{u}(t))(y_{}
^{}(t)+z(t)) \\
&+\frac{1}{2}p(t)^{\rm{T}}b_{xx}(\bar{X}(t),\bar{u}(t))(y(t))^2
+p(t)^{\rm{T}}(b(\bar{X}{(t)},u^{\varepsilon}(t))-b(\bar{X}{(t)},\bar{u}(t)))
\\
&+\displaystyle\sum_{j=1}^d\big{(}\frac{1}{2}q^j(t)^{\rm{T}}\sigma_{xx}^j(\bar{X}(t),\bar{u}(t))(y(t))^2
+q^j(t)^{\rm{T}}(\sigma(\bar{X}{(t)},u^{\varepsilon}(t))-\sigma(\bar{X}{(t)},\bar{u}(t)))\big{)}\bigg{]}dt.
\end{array}
\]
Therefore,
\begin{equation}%
\begin{array}
[c]{cl}
&
E\big{[}-\displaystyle \sum \limits_{i=1}^{n}
\phi_x(\bar{X}(t_i))(y^{}(t_{i})+z(t_i))\big{]}\\
= &
E{\displaystyle \int \limits_{0}^{T}}
\bigg{[}f_x(\bar{X}{(t)},\bar{u}(t))(y_{}
^{}(t)+z(t)) \\
&+\frac{1}{2}p(t)^{\rm{T}}b_{xx}(\bar{X}(t),\bar{u}(t))(y(t))^2
+p(t)^{\rm{T}}(b(\bar{X}{(t)},u^{\varepsilon}(t))-b(\bar{X}{(t)},\bar{u}(t)))
\\
&+\displaystyle\sum_{j=1}^d\big{(}\frac{1}{2}q^j(t)^{\rm{T}}\sigma_{xx}^j(\bar{X}(t),\bar{u}(t))(y(t))^2
+q^j(t)^{\rm{T}}(\sigma(\bar{X}{(t)},u^{\varepsilon}(t))-\sigma(\bar{X}{(t)},\bar{u}(t)))\big{)}\bigg{]}dt.
\end{array}
\label{maxm-2}%
\end{equation}
Now, let $u(t)=u$ be a constant, and $E_{\varepsilon}=[v,v+\varepsilon
]\subset \lbrack0,T]$. Combining equation (\ref{varm-3}) with (\ref{maxm-2}) and
noting the optimality of $\bar{u}(\cdot)$, we obtain%
\[%
\begin{array}
[c]{rl}%
0\leq & J(u^{\varepsilon}(\cdot))-J(\bar{u}(\cdot))\\
= &
\displaystyle \sum \limits_{i=1}^{n}E\left[\phi(\bar{X}(t_i))(y(t_{i})+z(t_{i}))
+\phi_{xx}(\bar{X}(t_i))y^{}(t_{i})y^{}(t_{i})\right]\\
& +
E{\displaystyle \int \limits_{0}^{T}}
\big{[}f_x(\bar{X}{(t)},\bar{u}(t))(y(t)+z(t))+\displaystyle\frac{1}{2}f_{xx}(\bar{X}{(t)},\bar{u}(t))(y(t))^2\\
&+f(\bar{X}{(t)},u^{\varepsilon}(t))-f(\bar{X}{(t)},\bar{u}(t))\big{]}dt+o(\varepsilon)
\\
=&E\displaystyle \sum \limits_{i=1}^{n}\phi_{xx}(\bar{X}(t_i))y^{}(t_{i})y^{}(t_{i})+
E{\displaystyle \int \limits_{0}^{T}}
\big{[}f_x(\bar{X}{(t)},\bar{u}(t))(y(t)+z(t))+\displaystyle \frac{1}{2}f_{xx}(\bar{X}{(t)},\bar{u}(t))(y(t))^2\\
&+f(\bar{X}{(t)},u^{\varepsilon}(t))-f(\bar{X}{(t)},\bar{u}(t))\big{]}dt+o(\varepsilon)\\
&-
E{\displaystyle \int \limits_{0}^{T}}
\bigg{[}f_x(\bar{X}{(t)},\bar{u}(t))(y_{}
^{}(t)+z(t)) \\
&+\displaystyle\frac{1}{2}p(t)^{\rm{T}}b_{xx}(\bar{X}(t),\bar{u}(t))(y(t))^2
+p(t)^{\rm{T}}(b(\bar{X}{(t)},u^{\varepsilon}(t))-b(\bar{X}{(t)},\bar{u}(t)))
\\
&+\displaystyle\sum_{j=1}^d\big{(}\frac{1}{2}q^j(t)^{\rm{T}}\sigma_{xx}^j(\bar{X}(t),\bar{u}(t))(y(t))^2
+q^j(t)^{\rm{T}}(\sigma(\bar{X}{(t)},u^{\varepsilon}(t))-\sigma(\bar{X}{(t)},\bar{u}(t)))\big{)}\bigg{]}dt\\
=&E\displaystyle \sum \limits_{i=1}^{n}\phi_{xx}(\bar{X}(t_i))y^{}(t_{i})y^{}(t_{i}) \\
&
+E\displaystyle\int_{0}^{T}\bigg{[}\frac{1}{2}f_{xx}(\bar{X}{(t)},\bar{u}(t))(y(t))^2
+f(\bar{X}{(t)},u^{\varepsilon}(t))-f(\bar{X}{(t)},\bar{u}(t))  \\
&-\displaystyle\{\frac{1}{2}p(t)^{\rm{T}}b_{xx}(\bar{X}(t),\bar{u}(t))(y(t))^2
+p(t)^{\rm{T}}(b(\bar{X}{(t)},u^{\varepsilon}(t))-b(\bar{X}{(t)},\bar{u}(t)))\\
&+\displaystyle\sum_{j=1}^n\big{(}\frac{1}{2}q^j(t)^{\rm{T}}\sigma_{xx}^j(\bar{X}(t),\bar{u}(t))(y(t))^2
+q^j(t)^{\rm{T}}(\sigma^j(\bar{X}{(t)},u^{\varepsilon}(t))-\sigma^j(\bar{X}{(t)},\bar{u}(t)))\big{)}\bigg{]}dt+o(\varepsilon).\\
\end{array}
\]
We recall that
\begin{equation*}
H(x,u,p,q)=b(x,u)^{\rm{T}}p+\displaystyle\sum_{j=1}^n\sigma^j (x,u)^{\rm{T}}q^j-f(x,u),\text{ \  \ }%
(x,u,p,q)\in \mathbb{R}^m\times U\times \mathbb{R}^m\times \mathbb{R}^{m\times d}.%
\end{equation*}
 Denote  $\bar{H}(t):=H(\bar{X}{(t)},\bar{u}(t),p(t),q(t))$ and $H^{\varepsilon}(t):=H(\bar{X}{(t)},{u}^{\varepsilon}(t),p(t),q(t))$.
  Thus,
\[%
\begin{array}
[c]{rl}%
 & J(u^{\varepsilon}(\cdot))-J(\bar{u}(\cdot))\\
= &E\displaystyle \sum \limits_{i=1}^{n}\phi_{xx}(\bar{X}(t_i))y^{}(t_{i})y^{}(t_{i}) -E{\displaystyle \int \limits_{0}^{T}}
\big{[}\displaystyle\frac{1}{2}\bar{H}_{xx}(t)(y(t))^2+H^{\varepsilon}(t)-\bar{H}(t)\big{]}dt+o(\varepsilon).\\
\end{array}
\]
Similar to the proof in Lemma \ref{lem-2}, by applying It\^{o} formula to $P(t)^{\rm{T}}(y(t))^2$ over $(t_{i-1},t_{i})$,  we have
\[%
\begin{array}
[c]{rl}%
 &E\big{[}\displaystyle \sum \limits_{i=1}^{n}\phi_{xx}(\bar{X}(t_i))y^{}(t_{i})y^{}(t_{i})\big{]} \\
= -&E {\displaystyle \int \limits_{0}^{T}}
\displaystyle\frac{1}{2}\big{[}\displaystyle\sum_{j=1}^d
(\sigma^j(\bar{X}{(t)},u^{\varepsilon}(t))-\sigma^j(\bar{X}{(t)},\bar{u}(t)))^{\rm{T}}
P(t)(\sigma(\bar{X}{(t)},u^{\varepsilon}(t))-\sigma(\bar{X}{(t)},\bar{u}(t)))\\
&+\bar{H}_{xx}(t)(y(t))^2\big{]}+o(\varepsilon).\\
\end{array}
\]
Then, we obtain
\[%
\begin{array}
[c]{rl}%
 0\geq &E {\displaystyle \int \limits_{0}^{T}} \big{[}H^{\varepsilon}(t)-\bar{H}(t)\\
 &+\displaystyle\frac{1}{2}\displaystyle\sum_{j=1}^d
(\sigma^j(\bar{X}{(t)},u^{\varepsilon}(t))-\sigma^j(\bar{X}{(t)},\bar{u}(t)))^{\rm{T}}
P(t)(\sigma(\bar{X}{(t)},u^{\varepsilon}(t))-\sigma(\bar{X}{(t)},\bar{u}(t))) \big{]}dt.\\
 \end{array}
\]
This completes the proof. $\ \ \ \ \ \ \ \ \Box$

\section{Maximum principle under  constraints}
We recall that
$$
\tau^n=\sum_{i=1}^{n-1}t_i1_{\{t_{i-1}\leq \tau < t_{i} \}}+t_n1_{\{t_{n-1}\leq \tau \leq  t_{n} \}},
$$
where $0=t_0< t_1<t_2\cdots<t_n=T$, and $t_i-t_{i-1}=\frac{T}{n},\ i=1,2,\cdots,n$. Similar to the proof in Lemma \ref{ner-le1}, for any given $\varepsilon>0$, there  exists a  large value of $n$, such that
$$
E\left|\Psi(X^u(\tau^n)) -\Psi(X^u(\tau))\right|<\varepsilon.
$$ All these conditions show that $E\Psi(X^u(\tau^n))$ is close to $E\Psi(X^u(\tau))$. Thus, we give the Pontryagin's  stochastic maximum principle  under  constrained conditions with stopping time $\tau^n$, i.e. the cost functional is as follows:
\begin{equation}
J(u(\cdot))=
E\displaystyle \int \limits_{0}^{T}
f(X{(t)},u(t))dt,\label{cost-3}%
\end{equation}
where the state process $X(\cdot)$ satisfies the following constrained condition:
\begin{equation}
\label{ccc-0}
0 \leq E\Psi(X^u(\tau^n)).
\end{equation}
Notice that
$$\Psi(X^u(\tau^n)=\sum_{i=1}^{n-1}\Psi(X^u(t_i))1_{\{t_{i-1}\leq \tau < t_{i} \}}+\Psi(X^u(t_n))1_{\{t_{n-1}\leq \tau \leq t_{n} \}}.$$
For notation simplicity, we set
\begin{equation}%
\begin{array}
[c]{rl}
&  \phi(X^u(t_i))=\Psi(X^u(t_i))1_{\{t_{i-1}\leq \tau <t_{i} \}},\ i=1,2,\cdots,n-1,\\
&  \phi(X^u(t_n))=\Psi(X^u(t_n))1_{\{t_{n-1}\leq \tau \leq t_{n} \}}.
\end{array}
\end{equation}
Thus, we can rewrite constrained condition (\ref{ccc-0}) as
\begin{equation}
\label{ccc-1}
0 \leq E\sum_{i=1}^n\phi(X^u(t_i)).
\end{equation}

To prove the main result of this section, we introduce the following Ekeland's variational principle, which comes from  Corollary 6.3 in \cite{YZ99}.
\begin{lemma}
\label{ccc-le1}
Let $F:V\to \mathbb{R}$ be a continuous function on a complete metric space $(V,\tilde{d})$. Given $\theta>0$ and $v_0\in V$ such that
$$
F(v_0)\leq \inf_{v\in V}F(v)+\theta.
$$
Then, there exists a $v_{\theta}\in V$ such that
$$
F(v_{\theta})\leq F(v_0),\ \ \tilde{d}(v_{\theta},v_0)\leq \sqrt{\theta},
$$
and for all $v\in V$,
$$
- \sqrt{\theta}d(v_{\theta},v)\leq F(v)-F(v_{\theta}).
$$
\end{lemma}

The related Hamiltonian is given as follows:
\begin{equation*}
H(\beta^0,x,u,p,q)=b(x,u)^{\text{T}}p+\sum_{j=1}^d\sigma^j (x,u)^{\text{T}}q^j-\beta^0f(x,u),\text{ \  \ }
\end{equation*}
where $(\beta^0,x,u,p,q)\in \mathbb{R}\times\mathbb{R}^m\times U\times \mathbb{R}^m\times \mathbb{R}^{m\times d}.$
Next, we present the main results of this section.

\begin{theorem}
\label{ccc-th}
Let Assumptions \ref{ass-b} and \ref{ass-fai} hold,
and $(\bar{u}(\cdot),\bar{X}(\cdot))$ be an optimal pair of (\ref{cost-3}).
Then, there exists $(\beta^0,\beta^1,\cdots,\beta^n)\in \mathbb{R}^{n+1}$ satisfying
$$
\beta^0\geq 0,\ \ \left| \beta^0\right|^2+\displaystyle\sum_{j=1}^n\left| \beta^j\right|^2=1,
$$
and
$$
\beta^j(\gamma-E\phi(\bar{X}(t_j)))\geq 0,\  \gamma\leq 0,\ j=1,2,\cdots,n.
$$
The adapted solution $(p(\cdot),q(\cdot))$ satisfies the following series of first-order adjoint equations:
\begin{equation}%
\begin{array}
[c]{ll}%
-d{p}(t)= & \big{[}b_x(\bar{X}{(t)},\bar{u}(t))^{\mathrm{T}}p(t)+ \sum_{j=1}^d\sigma_x^j(\bar{X}{(t)},\bar{u}(t))^{\text{T}}q^j(t) \\
               &-\beta^0f_x(\bar{X}{(t)},\bar{u}(t))\big{]}dt-q(t)dW(t),\ t\in(t_{i-1},t_{i}),\\
p(t_{i})= &-\beta^i\phi_{x}(\bar{X}(t_i)+p(t_{i}^{+}),\text{ \ }i=1,2,\ldots,n,
\end{array}
\label{prin-1}%
\end{equation}
and second-order adjoint equations,
\begin{equation}%
\begin{array}
[c]{ll}%
-dP(t)= & \bigg{\{}b_x(\bar{X}{(t)},\bar{u}(t))^{\rm{T}}P(t)+P(t)b_x(\bar{X}{(t)},\bar{u}(t))\\
              &+ \displaystyle\sum_{j=1}^d\sigma_x^j(\bar{X}{(t)},\bar{u}(t))^{\rm{T}} P(t)\sigma_x^j(\bar{X}{(t)},\bar{u}(t)) \\
               &+\displaystyle\sum_{j=1}^d\big{[}
               \sigma_x^j(\bar{X}{(t)},\bar{u}(t))^{\rm{T}}Q^j(t)+Q^j(t)\sigma_x^j(\bar{X}{(t)},\bar{u}(t))\big{]}\\
               & +H_{xx}(\bar{X}(t),\bar{u}(t),p(t),q(t))\bigg{\}}dt -Q(t)dW(t),\ t\in(t_{i-1},t_{i}),\\
-P(t_i)= & \beta^i\phi_{xx}(\bar{X}(t_i))-P(t_i^+),
\end{array}
\label{princ-2}%
\end{equation}
where $p(t_n^+)=0$, $P(t_n^+)=0$, and
\begin{equation}%
\begin{array}
[c]{ll}%
&H(\beta^0,\bar{X}(t),\bar{u}(t),p(t),q(t))-H(\beta^0,\bar{X}(t),u,p(t),q(t)))\\
\geq &\displaystyle\frac{1}{2} \displaystyle\sum_{j=1}^d\big{[}\sigma^j(\bar{X}{(t)},\bar{u}(t))-\sigma(\bar{X}{(t)},u)\big{]}^{\rm{T}}P(t)
 \big{[}\sigma^j(\bar{X}{(t)},\bar{u}(t))-\sigma(\bar{X}{(t)},u)\big{]},
\end{array}
\label{princ-3}%
\end{equation}
for any $u\in U$ and $t \in(t_{i},t_{i+1})$, where $i=0,1,\cdots,n-1$.
\end{theorem}
\noindent\textbf{Proof}: Without loss of generality, we assume that $J(\bar{u}(\cdot))=0$, where $(\bar{u}(\cdot),\bar{X}(\cdot))$ is the optimal pair of problem (\ref{cost-3}) with constrained condition (\ref{ccc-1}). For any $\theta>0$, we set
$$
J^{\theta}(u(\cdot))=\sqrt{\big{[}(J(u(\cdot))+\theta)^+\big{]}^2
+\displaystyle\sum_{i=1}^n\big{[}(-E\phi(X^u(t_i)))^+\big{]}^2}.
$$
In addition, one can verify that $J^{\theta}:\mathcal{U}[0,T]\to \mathbb{R}$ is a continuous function and satisfies
\begin{equation}
J^{\theta}(\bar{u}(\cdot))=\theta\leq \inf_{u\in\mathcal{U}[0,T]}J^{\theta}(u(\cdot))+\theta.
\end{equation}
Now, by Lemma \ref{ccc-le1}, there exists a $u^{\theta}(\cdot)\in \mathcal{U}[0,T]$ such that
\begin{equation}
\label{ccci-ev}
J^{\theta}(u^{\theta}(\cdot))\leq J^{\theta}(\bar{u}(\cdot))=\theta,\ \tilde{d}(u^{\theta}(\cdot),\bar{u}(\cdot))\leq \sqrt{\theta},
\end{equation}
where $\tilde{d}(u^1(\cdot),u^2(\cdot))=M\{(t,\omega)\in[0,T]\times \Omega:u^1(t,\omega)\neq u^2(t,\omega)\}$, $M$ is the product measure of the Lebesgue measure and probability on the set of $[0,T]\times \Omega$. We can verify that $(\mathcal{U}[0,T],\tilde{d})$ is a complete metric space. Also, we have
\begin{equation*}
-\sqrt{\theta}\tilde{d}(u^{\theta}(\cdot),u(\cdot))\leq J^{\theta}(u(\cdot))-J^{\theta}(u^{\theta}(\cdot)),\ \forall u(\cdot) \in \mathcal{U}[0,T],
\end{equation*}
which deduces
\begin{equation}
\label{cos-40}
J^{\theta}(u^{\theta}(\cdot))+\sqrt{\theta}\tilde{d}(u^{\theta}(\cdot),u^{\theta}(\cdot))\leq J^{\theta}(u(\cdot))+\sqrt{\theta}\tilde{d}(u^{\theta}(\cdot),u(\cdot)),\ \forall u(\cdot) \in \mathcal{U}[0,T].
\end{equation}
Thus, inequality (\ref{cos-40}) shows that $(u^{\theta}(\cdot),X^{\theta}(\cdot))$ is the optimal pair for the cost functional
\begin{equation}
\label{cos-4}
J^{\theta}(u(\cdot))+\sqrt{\theta}\tilde{d}(u^{\theta}(\cdot),u(\cdot)),
\end{equation}
without the state constraints.

Since $U$ is a general control domain, let $\rho>0$ and $E_{\rho}=[v,v+\rho
]\subset (t_{i-1},t_i)$,  for some $i\in\{1,2,\cdots,n\}$. Let $u\in U$ be any given constant. We define the following:%
\[
u^{\theta,\rho}(t)=\left \{
\begin{array}
[c]{cl}%
{u}^{\theta}(t), & \text{if }t\in \lbrack0,T]\backslash E_{\rho},\\
u, & \text{if }t\in E_{\rho},
\end{array}
\right.
\]
which belongs to $ \mathcal{U}[0,T]$. It is easy to verify that
$$
\tilde{d}(u^{\theta,\rho}(\cdot),u^{\theta}(\cdot))\leq \rho.
$$
By applying equation (\ref{cos-40}), we get
\begin{equation}%
\begin{array}
[c]{rl}%
-\sqrt{\theta}\rho \leq & J^{\theta}(u^{\theta,\rho}(\cdot))-J^{\theta}(u^{\theta}(\cdot))\\
=&\displaystyle\frac{\big{[}(J(u^{\theta,\rho}(\cdot))+\theta)^+\big{]}^2-\big{[}(J(u^{\theta}(\cdot))+\theta)^+\big{]}^2
}
{J^{\theta}(u^{\theta,\rho}(\cdot))+J^{\theta}(u^{\theta}(\cdot))}\\
&+\displaystyle\frac{\sum_{j=1}^n\big{[}\big{[}(-E\phi(X^{\theta,\rho}(t_j)))^+\big{]}^2-
\big{[}(-E\phi(X^{\theta}(t_j)))^+\big{]}^2\big{]}}
{J^{\theta}(u^{\theta,\rho}(\cdot))+J^{\theta}(u^{\theta}(\cdot))},\\
\end{array}
\label{ccci-1}%
\end{equation}
where $X^{\theta,\rho}(\cdot)$ and $X^{\theta}(\cdot)$ are the related solutions of equation (\ref{ODE_1}) with controls $u^{\theta,\rho}(\cdot)$ and $u^{\theta}(\cdot)$. We set
\begin{equation}%
\begin{array}
[c]{ll}%
\beta^{0,\theta}=\displaystyle\frac{\big{[}J(u^{\theta}(\cdot))+\theta\big{]}^+}{J^{\theta}(u^{\theta}(\cdot))},\\
\beta^{j,\theta}=\displaystyle\frac{-\big{[}-E\phi(X^{\theta}(t_j))\big{]}^+}
{J^{\theta}(u^{\theta}(\cdot))},\ j=1,2,\cdots,n.\\
\end{array}
\label{ccci-2}%
\end{equation}
Then, by the continuity of $J^{\theta}(\cdot)$, we have
\begin{equation}%
\begin{array}
[c]{rl}%
& J^{\theta}(u^{\theta,\rho}(\cdot))-J^{\theta}(u^{\theta}(\cdot))\\
=&\beta^{0,\theta}\big{[}J(u^{\theta,\rho}(\cdot))-J(u^{\theta}(\cdot))\big{]}+
\displaystyle\sum_{j=1}^n\beta^{j,\theta}\big{[}E\phi(X^{\theta,\rho}(t_j))
-E\phi(X^{\theta}(t_j))\big{]}+o(1),\\
=&E\big{[}\displaystyle\sum_{j=1}^{n}\beta^{j,\theta}
(\phi(X^{\theta,\rho}(t_j))
-\phi(X^{\theta}(t_j)))\\
&+\beta^{0,\theta}\displaystyle\int_0^T\big{[}f(X^{\theta,\rho}(t),u^{\theta,\rho}(t))-
f(X^{\theta}(t),u^{\theta}(t))\big{]}dt   \big{]}+o(\rho),
\end{array}
\label{ccci-3}%
\end{equation}
where $o(1)$ converges to $0$ when $\rho\to 0$.

Similar to the proof in Lemma \ref{lem-2}, let $(\bar{X}(\cdot),\bar{u}(\cdot))$ be replaced by $(X^{\theta}(t),u^{\theta}(t))$, $y(\cdot)$ be replaced by $\tilde{y}(\cdot)$ in equation (\ref{aprom-1}), and $z(\cdot)$ be replaced by $\tilde{z}(\cdot)$ in equation (\ref{aprom-2}). Thus, we obtain
\begin{equation}%
\begin{array}
[c]{rl}
-\sqrt{\theta} \rho\leq & J^{\theta}(u^{\theta,\rho}(\cdot))-J^{\theta}(u^{\theta}(\cdot))\\
\leq&E\bigg{[}
\displaystyle\sum_{i=1}^{n}\beta^{i,\theta}\phi_x(X^{\theta}(t_i))(\tilde{y}(t_i)+\tilde{z}(t_i))\\
&+\displaystyle\sum_{i=1}^{n}\beta^{i,\theta}\frac{1}{2}\phi_{xx}(X^{\theta}(t_i))(\tilde{y}(t_i),\tilde{y}(t_i))\\
&+\beta^{0,\theta}{\displaystyle \int \limits_{0}^{T}}
\big{(}f_x(X^{\theta}{(t)},{u}^{\theta}(t))(\tilde{y}(t)+\tilde{z}(t))
+\frac{1}{2}f_{xx}({X}^{\theta}{(t)},{u}^{\theta}(t))(\tilde{y}(t))^2\\
&+f(X^{\theta}{(t)},u^{\theta,\rho}(t))-f(X^{\theta}{(t)},{u}^{\theta}(t))\big{)}dt\bigg{]}+o(\rho).\\
\end{array}
\label{ccci-4}%
\end{equation}
In addition, we introduce the following adjoint equations:
\begin{equation}%
\begin{array}
[c]{rl}%
-d{p}^{\theta}(t)= & \big{[}b_x(X^{\theta}(t),u^{\theta}(t))^{\mathrm{T}}p^{\theta}(t)+ \sum_{j=1}^d\sigma_x^j(X^{\theta}(t),u^{\theta}(t))^{\text{T}}q^{j,{\theta}}(t) \\
               &-\beta^{0,\theta}f_x(X^{\theta}(t),u^{\theta}(t))\big{]}dt-q^{\theta}(t)dW(t),\ t\in(t_{i-1},t_{i}),\\
p^{\theta}(t_{i})= &-\beta^{i,\theta}\phi_x(X^{\theta}(t_i))+p(t_{i}^{+}),\text{ \ }i=1,\ldots,n,
\end{array}
\label{ccci-5}%
\end{equation}
where $q^{\theta}(\cdot)=(q^{1,\theta}(\cdot),q^{2,\theta}(\cdot),\cdots,q^{d,\theta}(\cdot))$, and
\begin{equation}%
\begin{array}
[c]{ll}%
-dP^{\theta}(t)= &\bigg{ \{}b_x(X^{\theta}(t),u^{\theta}(t))^{\rm{T}}P^{\theta}(t)
+P^{\theta}(t)b_x(X^{\theta}(t),u^{\theta}(t))\\
              &+ \displaystyle\sum_{j=1}^d\sigma_x^j(X^{\theta}(t),u^{\theta}(t))^{\rm{T}} P^{\theta}(t)\sigma_x^j(X^{\theta}(t),u^{\theta}(t)) \\
               &+\displaystyle\sum_{j=1}^d\big{[}
               \sigma_x^j(X^{\theta}(t),u^{\theta}(t))^{\rm{T}}Q^{j,\theta}(t)
               +Q^{j,\theta}(t)\sigma_x^j(X^{\theta}(t),u^{\theta}(t))\big{]}\\
               & +H_{xx}(\beta^{0,\theta},X^{\theta}(t),u^{\theta}(t),p^{\theta}(t),q^{\theta}(t))\bigg{\}}dt -Q^{\theta}(t)dW(t),\\
-P^{\theta}(t_i)= & \beta^{i,\theta}\phi_{xx}(\bar{X}(t_i))-P^{\theta}(t_i^+),
\end{array}
\label{princm-2}%
\end{equation}
where $P^{\theta}(t_n^+)=0$ and  $Q^{\theta}(\cdot)=(Q^{1,\theta}(\cdot),Q^{2,\theta}(\cdot),\cdots,Q^{d,\theta}(\cdot))$. Now, by applying the duality relation as in the proof of Theorem \ref{Maximumprinciple}, we get
\[%
\begin{array}
[c]{rl}%
 &o(1)+\sqrt{\theta}\\
 \geq &E {\displaystyle \int \limits_{0}^{T}} \big{[}H^{\theta,\rho}(t,{u}^{\theta,\rho}(t))-{H}^{\theta}(t,{X}^{\theta}{(t)})\\
 &+\displaystyle\frac{1}{2}\displaystyle\sum_{j=1}^d
(\sigma^j({X}^{\theta}{(t)},u^{\theta,\rho}(t))-
\sigma^j({X}^{\theta}{(t)},{u}^{\theta,\rho}(t)))^{\rm{T}}
P^{\theta}(t)(\sigma({X}^{\theta}{(t)},u^{\theta,\rho}(t))-\sigma({X}^{\theta}{(t)},{u}^{\theta}(t))) \big{]}dt,\\
 \end{array}
\]
where
$$
{H}^{\theta}(t,{u}^{\theta}{(t)})=H(\beta^{0,\theta},{X}^{\theta}{(t)},{u}^{\theta}(t),p^{\theta}(t),q^{\theta}(t)),
$$
and
$$
H^{\theta,\rho}(t,{u}^{\theta,\rho}(t))=H(\beta^{0,\theta},{X}^{\theta}{(t)},{u}^{\theta,\rho}(t),p^{\theta}(t),q^{\theta}(t)).
$$
Notice that $o(1) \to 0$ when $\rho\to 0$. Thus, when $\rho\to 0$, we get
\begin{equation}
\label{ccci-6}
\begin{array}
[c]{rl}%
 \sqrt{\theta}
 \geq &H^{\theta}(t,u)-{H}^{\theta}(t,u^{\theta}(t))\\
 &+\displaystyle\frac{1}{2}\displaystyle\sum_{j=1}^d
(\sigma^j({X}^{\theta}{(t)},u)-
\sigma^j({X}^{\theta}{(t)},{u}^{\theta,\rho}(t)))^{\rm{T}}
P^{\theta}(t)(\sigma({X}^{\theta}{(t)},u)-
\sigma({X}^{\theta}{(t)},{u}^{\theta}(t))).\\
 \end{array}
\end{equation}

From inequality (\ref{ccci-ev}), we observe that $u^{\theta}(\cdot)$ converges  to $\bar{u}(\cdot)$ under $\tilde{d}$ as $\theta\to 0$. Then, by Assumptions \ref{ass-b}, \ref{assb-b2}, and the basic theory of stochastic differential equation, we have
$$
\displaystyle \sup_{0\leq t\leq T}E\left|X^{\theta}(t)-\bar{X}(t)\right| \to 0,
$$
as $\theta\to 0$. From equation (\ref{ccci-2}), we have
\begin{equation}
\label{unq-1}
\left|\beta^{0,\theta}\right|^2+\displaystyle\sum_{j=1}^n\left|\beta^{j,\theta}\right|^2=1.
\end{equation}
Thus, we can choose a sequence $\{\theta_k\}_{k=1}^{\infty}$ satisfying $\displaystyle\lim_{k\to\infty}\theta_k=0$ such that the limitations of $\beta^{0,\theta_k}$ and $\beta^{j,\theta_k}$ exist and we denote them by
\begin{equation}
\begin{array}
[c]{ll}
\beta^{0}=\displaystyle\lim_{\to\infty}\beta^{0,\theta_k},\\
\beta^{j}=\displaystyle\lim_{\to\infty}\beta^{j,\theta_k},\\
\end{array}
\end{equation}
respectively, with $j=1,2,\cdots,n$. From equation (\ref{unq-1}), we have
$$
\left|\beta^{0}\right|^2+\displaystyle\sum_{j=1}^n\left|\beta^{j}\right|^2=1,
$$
and
$$
\beta^j(\gamma-E\phi(\bar{X}(t_j)))\geq 0,\  \gamma\geq 0,\ j=1,2,\cdots,n.
$$
Similarly, we can obtain
$$
\displaystyle \sup_{0\leq t\leq T}E\big{[}\left|p^{\theta_k}(t)-p(t)\right|^2+\int_0^T\left|q^{\theta_k}(t)-q(t)\right|^2  dt\big{]}\to 0,
$$
as $k\to\infty$. When $k\to \infty$, from equation (\ref{ccci-6}), we have
\begin{equation}%
\begin{array}
[c]{ll}%
&H(\beta^0,\bar{X}(t),\bar{u}(t),p(t),q(t))-H(\beta^0,\bar{X}(t),u,p(t),q(t)))\\
\geq &\displaystyle\frac{1}{2} \displaystyle\sum_{j=1}^d\big{[}\sigma^j(\bar{X}{(t)},\bar{u}(t))-\sigma(\bar{X}{(t)},u)\big{]}^{\rm{T}}P(t)
 \big{[}\sigma^j(\bar{X}{(t)},\bar{u}(t))-\sigma(\bar{X}{(t)},u)\big{]},
\end{array}
\end{equation}
for any $u\in U$ and $t \in(t_{i},t_{i+1})$, $i=0,1,\cdots,n-1$.

Thus, we complete this proof. $\ \ \ \ \ \ \ \ \ \ \ \  \Box$

Next, we provide an example to illustrate the optimal production planning problem  under stopping time state constraints.

\begin{example}
Let $T=1$, and consider the following controlled stochastic differential equation:
\begin{equation}
\label{3exp-1}
\begin{array}
[c]{ll}
{X}^u(s)=\displaystyle\int_0^s\big{[}u(t)-y(t)\big{]}dt,\\
\end{array}
\end{equation}
where $y(\cdot)$ denotes an uncertain demand
$$
y(s)=\displaystyle\frac{8}{3} s- W(s),
$$
and $u(\cdot)=\{u(s),0\leq s\leq 1\}$ is a control process taking values in a
compact set $U=[0,2]$.
Thus, we minimise the following cost functional:
\begin{equation}
J(u(\cdot))=E[X^u(1)], \label{1exa-1}%
\end{equation}
with the state constraints
\begin{equation}
\label{eexa-1}
0\leq EX^u(\tau),\ EX^u(1),\ \,0\leq  \tau\leq 0.5.
\end{equation}
For a given integer $N$, let $\tau=\sum_{i=1}^{N-1}\frac{i}{2N}1_{\frac{i-1}{2N}\leq i< \frac{i}{2N}}+\frac{1}{2}1_{\frac{N-1}{2N}\leq i\leq \frac{1}{2}}$. Substituting  $X^u(\cdot)$ into equation (\ref{1exa-1}), we obtain
$$
J(u(\cdot))=E[
\displaystyle\int_{0}^1(u(t)-\frac{8}{3}t)dt].
$$
It is easy to verify that
 \begin{equation}
   (\tilde{u}^N(t),\tilde{X}^N(t))=
   \begin{cases}
    (\frac{2i-1}{4N},\displaystyle \int_0^tW(s)ds) &\mbox{ $\frac{i-1}{2N}\leq t \leq \frac{i}{2N}$,}\\
  (2,2t-\frac{4}{3}t^2-\frac{2}{3}+\displaystyle \int_0^tW(s)ds)&\mbox{ $ \frac{1}{2}  < t \leq 1$, }
   \end{cases}
  \end{equation}
is  an optimal pair of system (\ref{1exa-1}) under state constraints (\ref{eexa-1}).

Next, we show the maximum principle for the optimal control under multi-time state constraints (\ref{eexa-1}). Notice that it is difficult to obtain the adjoint equations for state process (\ref{3exp-1}) directly. To get the related adjoint equations, we rewrite equation (\ref{3exp-1}) as follows:
$$
X^u(s)-W(s)s=\displaystyle \int_0^s(u(t)-\frac{8}{3}t)dt-\displaystyle\int_0^stdW(t).
$$
To denote the above equation by $\delta X^u(s)=X^u(s)-W(s)s$, we have
$$
d\delta X^u(s)=[u(s)-\frac{8}{3}s]ds-sdW(s),
$$
and
 $$
 E[\delta X^u(1)]=E[X^u(1)],
 $$
which means that $(\bar{u}(t),\bar{X}(t)-W(s)s)$ is the optimal pair of the following cost functional:
\begin{equation}
\label{3exc-2}
\delta J(u(\cdot))=E[\delta X^u(1)],
\end{equation}
with the following constrained conditions:
$$
0\leq  E[\delta X^u(t_i)],E[\delta X^u(1)], \ t_i=\frac{i}{2N},\ i=1,2,\cdots,N.
$$
Next, we introduce the following  first-order adjoint  equations for functional (\ref{3exc-2}):
\begin{equation*}%
\begin{array}
[c]{ll}%
d{p}(t)= q(t)dW(t),\quad \frac{1}{2}  < t < 1,\\
p(1)= -(\beta^0+\beta^{N+1}),
\end{array}
\label{prin-11}%
\end{equation*}
 and
 \begin{equation*}%
 \begin{array}
[c]{ll}%
d{p}(t)= q(t)dW(t),\quad  \frac{i-1}{2N}< t < \frac{i}{2N}, \\
p(t_i)= -\beta^i+p(t_i^{+}),
\end{array}
\label{prin-12}%
\end{equation*}
and second-order adjoint equations
\begin{equation}%
\begin{array}
[c]{ll}%
dP(t)= & Q(t)dW(t),\\
P(T)= & 0,
\end{array}
\label{eprinc-2}%
\end{equation}
where $(\beta^0,\beta^1,\cdots,\beta^{N+1})$ is obtained from Theorem \ref{ccc-th}. The solutions of first- and second-order adjoint equations are given as follows:
 \begin{equation}
   (p(t),q(t))=
   \begin{cases}
   (-(\beta^0+\displaystyle\sum_{j=i}^{N+1}\beta^i),0) &\mbox{ $ \frac{i-1}{2N}< t \leq  \frac{i}{2N}$, }\\
  (-(\beta^0+\beta^{N+1}),0)&\mbox{ $\frac{1}{2} <  t \leq 1$,}
   \end{cases}
  \end{equation}
  and
 \begin{equation}
   (P(t),Q(t))=(0,0).
  \end{equation}
 Now, from Theorem \ref{ccc-th}, we have   $\beta^0=-\frac{\sqrt{2}}{2},\beta^{N+1}=\frac{\sqrt{2}}{2},\beta^{i}=0,i=1,2,\cdots,N$. From $\beta^0+\beta^{N+1}\leq0$ and $\beta^0+\displaystyle\sum_{j=i}^{N+1}\beta^i= 0,i=1,2,\cdots,N$, we obtain
 \begin{equation}%
\begin{array}
[c]{ll}%
&H(\beta^0,\bar{X}(t),\bar{u}(t),p(t),q(t))-H(\beta^0,\bar{X}(t),u,p(t),q(t)))\\
&-\displaystyle\frac{1}{2} \displaystyle\sum_{j=1}^d\big{(}\sigma^j(\bar{X}{(t)},\bar{u}(t))-\sigma(\bar{X}{(t)},u)\big{)}^{\rm{T}}P(t)
 \big{(}\sigma^j(\bar{X}{(t)},\bar{u}(t))-\sigma(\bar{X}{(t)},u)\big{)}\\
 =&(\tilde{u}^N(t)-u)p(t).
\end{array}
\label{eprinc-3}%
\end{equation}
Finally,
\begin{equation}
   (\tilde{u}^N(t)-u)p(t)=
   \begin{cases}
    -(\beta^0+\displaystyle\sum_{j=i}^{N+1}\beta^i))(\frac{2i-1}{4N}-u) &\mbox{ $ \frac{i-1}{2N}< t \leq  \frac{i}{2N}$, }\\
  -(\beta^0+\beta^{N+1}) (2-u) &\mbox{ $\frac{1}{2} < t \leq 1$,}\\
\end{cases}
     \end{equation}
 where $u\in[0,2]$.  Hence, the optimal control pair $(\tilde{u}^N(\cdot),\tilde{X}^N(\cdot))$ satisfies Theorem \ref{ccc-th}.
\end{example}

\section{Conclusion}
In this study, we considered the stochastic optimal control problem with a stopping time cost functional. We considered the classical constant terminal time as a stopping time in the cost functional. For example, in the financial market, the time when the investor leaves the market can be considered as the stopping time. To solve this kind of problem, we introduced a near-optimal control problem for the stopping time cost functional. Based on a series of first- and second-order adjoint equations, we established a stochastic maximum principle for the near-optimal control problem and the near-optimal control problem under a stopping time state constraints.  A closely related work \cite{Y16c} provides the necessary and sufficient conditions for stochastic differential systems under multi-time states cost functional with a convex control domain. In future, certain related topics such as, to show the dynamic programming principle for this system, to solve the mean-variance problem, and to perform empirical analysis, should be considered for the stopping time cost functional.

\end{document}